\documentclass[conference]{IEEEtran}

\usepackage{cite}
\usepackage{graphicx}
\usepackage[breaklinks=false, allcolors=black, colorlinks=true]{hyperref}

\usepackage{amsfonts}
\usepackage{amsmath} 
\usepackage{amssymb}
\usepackage{amsthm}	
\usepackage{enumerate}
\usepackage{mathrsfs}
\usepackage{subcaption}
\usepackage{epstopdf}

\newtheorem{lemm}{Lemma}
\newtheorem{definitio}{Definition}

\newtheorem{propositio}{Proposition}

\newtheorem{remar}{Remark}

\def\nn{\nonumber}
\def\beq{\begin{equation}}
\def\eeq{\end{equation}}
\def\beqs{\begin{equation*}}
\def\eeqs{\end{equation*}}
\def\balign{\begin{align}}
\def\ealign{\end{align}}
\def\bea{\begin{eqnarray}}
\def\eea{\end{eqnarray}}
\def\ba{\begin{array}}
\def\ea{\end{array}}
\def\bcenter{\begin{center}}
\def\ecenter{\end{center}}
\def\bitem{\begin{itemize}}
\def\eitem{\end{itemize}}
\def\benum{\begin{enumerate}}
\def\eenum{\end{enumerate}}
\def\bdoc{\begin{document}}
\def\edoc{\end{document}}
\def\defeq{{\stackrel{\Delta}{=}}}
\def\nn{\nonumber}

\def\eg{{e.g., \/}}
\def\etal{{\it et al. \/}}
\def\viz{{\iet viz.,\ \/}}
\def\ie{{i.e.,\ \/}}
\def\etc{{\it etc.,\ \/}}
\def\cf{{cf.,\ \/}}

\usepackage[usenames,dvipsnames,svgnames,table]{xcolor}

\def\yellow{\textcolor{yellow}}
\def\green{\textcolor{green}}
\def\tcbg{\textcolor{bgrd}}
\def\magenta{\textcolor{magenta}}
\def\white{\textcolor{white}}
\def\gray{\textcolor{gray}}

\newcommand\Quote[1]{\lq\textsl{#1}\rq}
\newcommand\fr[2]{{\textstyle\frac{#1}{#2}}}

\newcommand{\Ambb}{\mathbb{A}}
\newcommand{\Bmbb}{\mathbb{B}}
\newcommand{\Cmbb}{\mathbb{C}}
\newcommand{\Dmbb}{\mathbb{D}}
\newcommand{\Embb}{\mathbb{E}}
\newcommand{\Fmbb}{\mathbb{F}}
\newcommand{\Gmbb}{\mathbb{G}}
\newcommand{\Hmbb}{\mathbb{H}}
\newcommand{\Imbb}{\mathbb{I}}
\newcommand{\Jmbb}{\mathbb{J}}
\newcommand{\Kmbb}{\mathbb{K}}
\newcommand{\Lmbb}{\mathbb{L}}
\newcommand{\Mmbb}{\mathbb{M}}
\newcommand{\Nmbb}{\mathbb{N}}
\newcommand{\Ombb}{\mathbb{O}}
\newcommand{\Pmbb}{\mathbb{P}}
\newcommand{\Qmbb}{\mathbb{Q}}
\newcommand{\Rmbb}{\mathbb{R}}
\newcommand{\Smbb}{\mathbb{S}}
\newcommand{\Tmbb}{\mathbb{T}}
\newcommand{\Umbb}{\mathbb{U}}
\newcommand{\Vmbb}{\mathbb{V}}
\newcommand{\Wmbb}{\mathbb{W}}
\newcommand{\Xmbb}{\mathbb{X}}
\newcommand{\Ymbb}{\mathbb{Y}}
\newcommand{\Zmbb}{\mathbb{Z}}

\newcommand{\Amsc}{\mathcal{A}}
\newcommand{\Bmsc}{\mathcal{B}}
\newcommand{\Cmsc}{\mathcal{C}}
\newcommand{\Dmsc}{\mathcal{D}}
\newcommand{\Emsc}{\mathcal{E}}
\newcommand{\Fmsc}{\mathcal{F}}
\newcommand{\Gmsc}{\mathcal{G}}
\newcommand{\Hmsc}{\mathcal{H}}
\newcommand{\Imsc}{\mathcal{I}}
\newcommand{\Jmsc}{\mathcal{J}}
\newcommand{\Kmsc}{\mathcal{K}}
\newcommand{\Lmsc}{\mathcal{L}}
\newcommand{\Mmsc}{\mathcal{M}}
\newcommand{\Nmsc}{\mathcal{N}}
\newcommand{\Omsc}{\mathcal{O}}
\newcommand{\Pmsc}{\mathcal{P}}
\newcommand{\Qmsc}{\mathcal{Q}}
\newcommand{\Rmsc}{\mathcal{R}}
\newcommand{\Smsc}{\mathcal{S}}
\newcommand{\Tmsc}{\mathcal{T}}
\newcommand{\Umsc}{\mathcal{U}}
\newcommand{\Vmsc}{\mathcal{V}}
\newcommand{\Wmsc}{\mathcal{W}}
\newcommand{\Xmsc}{\mathcal{X}}
\newcommand{\Ymsc}{\mathcal{Y}}
\newcommand{\Zmsc}{\mathcal{Z}}

\newcommand{\Ascr}{\mathscr{A}}
\newcommand{\Bscr}{\mathscr{B}}
\newcommand{\Cscr}{\mathscr{C}}
\newcommand{\Dscr}{\mathscr{D}}
\newcommand{\Escr}{\mathscr{E}}
\newcommand{\Fscr}{\mathscr{F}}
\newcommand{\Gscr}{\mathscr{G}}
\newcommand{\Hscr}{\mathscr{H}}
\newcommand{\Iscr}{\mathscr{I}}
\newcommand{\Jscr}{\mathscr{J}}
\newcommand{\Kscr}{\mathscr{K}}
\newcommand{\Lscr}{\mathscr{L}}
\newcommand{\Mscr}{\mathscr{M}}
\newcommand{\Nscr}{\mathscr{N}}
\newcommand{\Oscr}{\mathscr{O}}
\newcommand{\Pscr}{\mathscr{P}}
\newcommand{\Qscr}{\mathscr{Q}}
\newcommand{\Rscr}{\mathscr{R}}
\newcommand{\Sscr}{\mathscr{S}}
\newcommand{\Tscr}{\mathscr{T}}
\newcommand{\Uscr}{\mathscr{U}}
\newcommand{\Vscr}{\mathscr{V}}
\newcommand{\Wscr}{\mathscr{W}}
\newcommand{\Xscr}{\mathscr{X}}
\newcommand{\Yscr}{\mathscr{Y}}
\newcommand{\Zscr}{\mathscr{Z}}

\def\Atiny{\mbox{\tiny{A}}}
\def\Btiny{\mbox{\tiny{B}}}
\def\Ctiny{\mbox{\tiny{C}}}
\def\Dtiny{\mbox{\tiny{D}}}
\def\Etiny{\mbox{\tiny{E}}}
\def\Ctiny{\mbox{\tiny{C}}}
\def\Ftiny{\mbox{\tiny{F}}}
\def\Gtiny{\mbox{\tiny{G}}}
\def\Htiny{\mbox{\tiny{H}}}
\def\Itiny{\mbox{\tiny{I}}}
\def\Jtiny{\mbox{\tiny{J}}}
\def\Ktiny{\mbox{\tiny{K}}}
\def\Ltiny{\mbox{\tiny{L}}}
\def\Mtiny{\mbox{\tiny{M}}}
\def\Ntiny{\mbox{\tiny{N}}}
\def\Otiny{\mbox{\tiny{O}}}
\def\Ptiny{\mbox{\tiny{P}}}
\def\Qtiny{\mbox{\tiny{Q}}}
\def\Rtiny{\mbox{\tiny{R}}}
\def\Stiny{\mbox{\tiny{S}}}
\def\Ttiny{\mbox{\tiny{T}}}
\def\Utiny{\mbox{\tiny{U}}}
\def\Vtiny{\mbox{\tiny{V}}}
\def\Wtiny{\mbox{\tiny{W}}}
\def\Xtiny{\mbox{\tiny{X}}}
\def\Ytiny{\mbox{\tiny{Y}}}
\def\Ztiny{\mbox{\tiny{Z}}}

\def\alphabf{\hbox{\boldmath$\alpha$\unboldmath}}
\def\betabf{\hbox{\boldmath$\beta$\unboldmath}}
\def\gammabf{\hbox{\boldmath$\gamma$\unboldmath}}
\def\deltabf{\hbox{\boldmath$\delta$\unboldmath}}
\def\epsilonbf{\hbox{\boldmath$\epsilon$\unboldmath}}
\def\zetabf{\hbox{\boldmath$\zeta$\unboldmath}}
\def\etabf{\hbox{\boldmath$\eta$\unboldmath}}
\def\iotabf{\hbox{\boldmath$\iota$\unboldmath}}
\def\kappabf{\hbox{\boldmath$\kappa$\unboldmath}}
\def\lambdabf{\hbox{\boldmath$\lambda$\unboldmath}}
\def\mubf{\hbox{\boldmath$\mu$\unboldmath}}
\def\nubf{\hbox{\boldmath$\nu$\unboldmath}}
\def\xibf{\hbox{\boldmath$\xi$\unboldmath}}
\def\pibf{\hbox{\boldmath$\pi$\unboldmath}}
\def\rhobf{\hbox{\boldmath$\rho$\unboldmath}}
\def\sigmabf{\hbox{\boldmath$\sigma$\unboldmath}}
\def\taubf{\hbox{\boldmath$\tau$\unboldmath}}
\def\upsilonbf{\hbox{\boldmath$\upsilon$\unboldmath}}
\def\phibf{\hbox{\boldmath$\phi$\unboldmath}}
\def\chibf{\hbox{\boldmath$\chi$\unboldmath}}
\def\psibf{\hbox{\boldmath$\psi$\unboldmath}}
\def\omegabf{\hbox{\boldmath$\omega$\unboldmath}}
\def\Sigmabf{\hbox{$\bf \Sigma$}}
\def\Upsilonbf{\hbox{$\bf \Upsilon$}}
\def\Omegabf{\hbox{$\bf \Omega$}}
\def\Deltabf{\hbox{$\bf \Delta$}}
\def\Gammabf{\hbox{$\bf \Gamma$}}
\def\Thetabf{\hbox{$\bf \Theta$}}
\def\Lambdabf{\mbox{$ \bf \Lambda $}}
\def\Xibf{\hbox{\bf$\Xi$}}
\def\Pibf{{\bf \Pi}}
\def\thetabf{{\mbox{\boldmath$\theta$\unboldmath}}}
\def\Upsilonbf{\hbox{\boldmath$\Upsilon$\unboldmath}}
\newcommand{\Phibf}{\mbox{${\bf \Phi}$}}
\newcommand{\Psibf}{\mbox{${\bf \Psi}$}}

\def\abf{{\bf a}}
\def\bbf{{\bf b}}
\def\cbf{{\bf c}}
\def\dbf{{\bf d}}
\def\ebf{{\bf e}}
\def\fbf{{\bf f}}
\def\gbf{{\bf g}}
\def\hbf{{\bf h}}
\def\ibf{{\bf i}}
\def\jbf{{\bf j}}
\def\kbf{{\bf k}}
\def\lbf{{\bf l}}
\def\mbf{{\bf m}}
\def\nbf{{\bf n}}
\def\obf{{\bf o}}
\def\pbf{{\bf p}}
\def\qbf{{\bf q}}
\def\rbf{{\bf r}}
\def\sbf{{\bf s}}
\def\tbf{{\bf t}}
\def\ubf{{\bf u}}
\def\vbf{{\bf v}}
\def\wbf{{\bf w}}
\def\xbf{{\bf x}}
\def\ybf{{\bf y}}
\def\zbf{{\bf z}}
\def\rbf{{\bf r}}
\def\xbf{{\bf x}}
\def\ybf{{\bf y}}
\def\Abf{{\bf A}}
\def\Bbf{{\bf B}}
\def\Cbf{{\bf C}}
\def\Dbf{{\bf D}}
\def\Ebf{{\bf E}}
\def\Fbf{{\bf F}}
\def\Gbf{{\bf G}}
\def\Hbf{{\bf H}}
\def\Ibf{{\bf I}}
\def\Jbf{{\bf J}}
\def\Kbf{{\bf K}}
\def\Lbf{{\bf L}}
\def\Mbf{{\bf M}}
\def\Nbf{{\bf N}}
\def\Obf{{\bf O}}
\def\Pbf{{\bf P}}
\def\Qbf{{\bf Q}}
\def\Rbf{{\bf R}}
\def\Sbf{{\bf S}}
\def\Tbf{{\bf T}}
\def\Ubf{{\bf U}}
\def\Vbf{{\bf V}}
\def\Wbf{{\bf W}}
\def\Xbf{{\bf X}}
\def\Ybf{{\bf Y}}
\def\Zbf{{\bf Z}}
\def\Ac{{\cal A}}
\def\Bc{{\cal B}}
\def\Cc{{\cal C}}
\def\Dc{{\cal D}}
\def\Ec{{\cal E}}
\def\Fc{{\cal F}}
\def\Gc{{\cal G}}
\def\Hc{{\cal H}}
\def\Ic{{\cal I}}
\def\Jc{{\cal J}}
\def\Kc{{\cal K}}
\def\Lc{{\cal L}}
\def\Mc{{\cal M}}
\def\Nc{{\cal N}}
\def\Oc{{\cal O}}
\def\Pc{{\cal P}}
\def\Qc{{\cal Q}}
\def\Rc{{\cal R}}
\def\Sc{{\cal S}}
\def\Tc{{\cal T}}
\def\Uc{{\cal U}}
\def\Vc{{\cal V}}
\def\Wc{{\cal W}}
\def\Xc{{\cal X}}
\def\Yc{{\cal Y}}
\def\Zc{{\cal Z}}

\def\0bf{\mathbf{0}}
\def\1bf{\mathbf{1}}

\def\lambdabf{\boldsymbol{\lambda}}
\def\mubf{\boldsymbol{\mu}}
\def\gammabf{\boldsymbol{\gamma}}
\def\xibf{\boldsymbol{\xi}}
\def\sigmabf{\boldsymbol{\sigma}}
\def\zetabf{\boldsymbol{\zeta}}

\def\Lambdabf{\boldsymbol{\Lambda}}
\def\Mubf{\boldsymbol{\Mu}}
\def\Gammabf{\boldsymbol{\Gamma}}
\def\Xibf{\boldsymbol{\Xi}}
\def\Sigmabf{\boldsymbol{\Sigma}}
\def\Zetabf{\boldsymbol{\Zeta}}

\usepackage{tikz}
\usetikzlibrary{arrows}

\newcommand\independent{\protect\mathpalette{\protect\independenT}{\perp}}
\def\independenT#1#2{\mathrel{\rlap{$#1#2$}\mkern2mu{#1#2}}}

\newcommand{\bsf}[1]{\textsf{\textbf{#1}}}
\newcommand{\lbsf}[1]{\textsf{\large  \textbf{#1}}}
\newcommand{\Lbsf}[1]{\textsf{\Large  \textbf{#1}}}
\newcommand{\hbsf}[1]{\textsf{\huge  \textbf{#1}}}

\newcommand{\myminipage}[3]{\begin{minipage}[#1]{#2}{#3} \end{minipage}}
\newcommand{\sbs}[4]{\myminipage{c}{#1}{#3} \hfill
\myminipage{c}{#2}{#4}}

\newcommand{\myfig}[2]{\centerline{\psfig{figure=#1,width=#2,silent=}}}
\newcommand{\myfigh}[2]{\centerline{\psfig{figure=#1,height=#2,silent=}}}
\newcommand{\myfigwh}[3]{\centerline{\psfig{figure=#1,width=#2,height=#3,silent=}}}

\newcommand{\beqa}{\begin{eqnarray}}
\newcommand{\eeqa}{\end{eqnarray}}
\newcommand{\beqan}{\begin{eqnarray*}}
\newcommand{\eeqan}{\end{eqnarray*}}
\newcommand{\dst}[1]{\displaystyle{ #1 }}
\newcounter{pf}

\newcommand{\smax}[1] { \bar \sigma \left( #1 \right) }
\newcommand{\Rn}{{\mathbb R}^n}
\newcommand{\R}{{\mathbb R}}
\newcommand{\C}{{\mathbb C}}
\newcommand{\Rm}{\mathbb{R}^m}
\newcommand{\Rmn}{\mathbb{R}^{m \times n}}
\newcommand{\Rpq}{\mathbb{R}^{p \times q}}
\newcommand{\Cn}{\mathbb{C}^n}
\newcommand{\Cm}{\mathbb{C}^m}
\newcommand{\Cnn}{\mathbb{C}^{n \times n}}
\newcommand{\Cmn}{\mathbb{C}^{m \times n}}
\newcommand{\ip}[1]{\left\langle #1 \right\rangle}
\newcommand{\rank}{\mbox{rank}}
\newcommand{\Span}{\mbox{\rm Span }}
\newcommand{\Trace}{\mbox{\rm Tr }}
\newcommand{\trace}[1]{\text{Tr}\left(#1\right)}
\newcommand{\Spec}{\mbox{\rm Spec }}
\newcommand{\vectornorm}[1]{\left\|#1\right\|}

\newcommand{\pd}[2]{\frac{\partial #1}{\partial #2}}
\newcommand{\ppd}[3]{\frac{\partial^2 #1}{\partial #2 \partial #3}}

\newcommand{\thtilde}{\tilde{\theta}}
\newcommand{\thnom}{\theta^\circ}
\newcommand{\thopt}{\theta^{\mbox{\small opt}}}
\newcommand{\thhat}{{\hat{\theta}}}
\newcommand{\Tho}{\Theta^\circ}
\newcommand{\tho}{\theta^\circ}
\newcommand{\np}{{n_p}}

\newcommand{\ii}{{[i]}}
\newcommand{\II}{{[i+1]}}
\newcommand{\iii}{{[ii]}}
\newcommand{\jj}{{[j]}}
\newcommand{\kk}{{[k]}}
\newcommand{\thi}{{\theta^\ii}}
\newcommand{\thI}{{\theta^\II}}
\newcommand{\di}{{d^\ii}}
\newcommand{\gi}{{g^\ii}}
\newcommand{\Hi}{{\HH^\ii}}
\newcommand{\thK}{\theta^{(k+1)}}
\newcommand{\gk}{{g^{(k)}}}
\newcommand{\Hk}{{{\cal H}^{(k)}}}

\newcommand{\bfdelta}{{\bf \Delta}}

\newcommand{\Exp}[1]{\exp \left\{ #1 \right\}} 
\newcommand{\gaussian}[1]{\mathbb{N} \left( #1 \right)}
\newcommand{\uniform}[1]{\mathbb{U} \left[ #1 \right]}
\newcommand{\exponential}[1]{\mathbb{E} \left[ #1 \right]}
\newcommand{\EXP}[1]{\EEXP \left[ #1 \right]} 
\newcommand{\EEXP}{\mbox{\bsf{E}}} 
\newcommand{\Prob}[1]{\mbox{{\sf Pr}} \left(#1 \right)}
\newcommand{\convas}{\stackrel{as}{\longrightarrow}}
\newcommand{\convinp}{\stackrel{p}{\longrightarrow}}
\newcommand{\convind}{\stackrel{d}{\longrightarrow}}
\newcommand{\convqm}{\stackrel{qm}{\longrightarrow}}
\newcommand{\sss}[1]{{_{#1}}}
\newcommand{\density}[2]{p_{_{_{#1}}}\!\!\left(#2 \right)} 
\newcommand{\distro}[2]{P_{_{_{#1}}}\!\!\left(#2 \right)} 
\newcommand{\rxx}[1]{R_{_{#1}}\!} 
\newcommand{\sxx}[1]{S_{_{#1}}} 
\newcommand{\cov}[1]{\Lambda_{_{#1}}} 
\newcommand{\mean}[1]{m_{_{#1}}} 
\newcommand{\LS}[1]{\hat{#1}_{_{LS}}} 
\newcommand{\MV}[1]{\hat{#1}_{_{MV}}} 
\newcommand{\LMV}[1]{\hat{#1}_{_{LMV}}} 
\newcommand{\ML}[1]{\hat{#1}_{_{ML}}} 

\renewcommand{\arraystretch}{0.9}
\newcommand{\bmat}[1]{ \begin{bmatrix} #1 \end{bmatrix}}
\newcommand{\mat}[1]{ \left[ \begin{array}{cccccccc} #1 \end{array}
\right] }
\newcommand{\smallmat}[1]{\small{\mat{#1}}}
\newcommand{\sysblk}[4]{\begin{array}{c|cccc}#1&#2\\ \hline#3&#4
\end{array}}
\newcommand{\sysmat}[4]{\left[\sysblk{#1}{#2}{#3}{#4}\right]}
\newcommand{\SGeq}{\succ}
\newcommand{\SLeq}{\prec}
\newcommand{\Geq}{\succeq}
\newcommand{\Leq}{\preceq}

\newcommand{\Aset}{\mathbb{A}}
\newcommand{\Bset}{\mathbb{B}}
\newcommand{\Fset}{\mathbb{F}}
\newcommand{\Gset}{\mathbb{G}}
\newcommand{\Kset}{\mathbb{K}}
\newcommand{\Mset}{\mathbb{M}}
\newcommand{\Sset}{\mathbb{S}}
\newcommand{\Tset}{\mathbb{T}}
\newcommand{\Uset}{\mathbb{U}}
\newcommand{\Vset}{\mathbb{V}}
\newcommand{\Wset}{\mathbb{W}}

\newcommand{\Ical}{{\cal I}}
\newcommand{\Acal}{{\cal A}}
\newcommand{\Bcal}{{\cal B}}
\newcommand{\Ccal}{{\mathcal{C}}}
\newcommand{\Dcal}{{\cal D}}
\newcommand{\Ecal}{{\mathcal{E}}}
\newcommand{\Fcal}{{\cal F}}
\newcommand{\Gcal}{{\mathcal{G}}}
\newcommand{\Hcal}{{\cal H}}
\newcommand{\Kcal}{{\mathcal{K}}}
\newcommand{\Lcal}{{\cal L}}
\newcommand{\Mcal}{{\cal M}}
\newcommand{\Ncal}{{\mathcal{N}}}
\newcommand{\Pcal}{{\cal P}}
\newcommand{\Qcal}{{\mathcal{Q}}}
\newcommand{\Rcal}{{\cal R}}
\newcommand{\Scal}{{\mathcal{S}}}
\newcommand{\Tcal}{{\mathcal{T}}}
\newcommand{\Wcal}{{\mathcal{W}}}
\newcommand{\Ucal}{{\cal U}}
\newcommand{\Vcal}{{\mathcal{V}}}
\newcommand{\Xcal}{{\cal X}}
\newcommand{\Ycal}{{\mathcal{Y}}}
\newcommand{\Zcal}{{\mathcal{Z}}}

\newcommand{\EE}{{\bf E}}
\newcommand{\FF}{{\bf F}}
\newcommand{\GG}{{\bf G}}
\newcommand{\HH}{{\bf H}}
\newcommand{\LL}{{\bf L}}
\newcommand{\NN}{{\bf N}}
\newcommand{\MM}{{\bf M}}
\newcommand{\PP}{{\bf P}}
\newcommand{\QQ}{{\bf Q}}
\newcommand{\RR}{{\bf R}}
\renewcommand{\SS}{{\bf S}}
\newcommand{\TT}{{\bf T}}
\newcommand{\VV}{{\bf V}}
\newcommand{\WW}{{\bf W}}

\newcommand{\thk}{\theta^{(k)}}
\newcommand{\thb}{\theta^{\rm opt}}
\newcommand{\alb}{\alpha^{\rm opt}}
\newcommand{\dk}{d^{(k)}}
\newcommand{\Hinf}{{\cal H}_\infty}
\newcommand{\Htwo}{{\cal H}_2}

\renewcommand{\arraystretch}{1.1}

\newcommand{\red}[1]{{\color{red} #1}}
\newcommand{\blue}[1]{{\color{Blue} #1}}
\newcommand{\black}[1]{{\color{Black} #1}}

\newcounter{l1}
\newcounter{l2}
\newcounter{l3}
\newcommand{\bdotlist}{\begin{list}{$\bullet$}{}}
\newcommand{\bboxlist}{\begin{list}{$\Box$}{}}
\newcommand{\bbboxlist}{\begin{list}{\raisebox{.005in}{{\tiny
$\blacksquare$ \ \ }}}{}}
\newcommand{\bdashlist}{\begin{list}{$-$}{} }
\newcommand{\blist}{\begin{list}{}{} }
\newcommand{\barablist}{\begin{list}{\arabic{l1}}{\usecounter{l1}}}
\newcommand{\balphlist}{\begin{list}{(\alph{l2})}{\usecounter{l2}}}
\newcommand{\bAlphlist}{\begin{list}{\Alph{l2}.}{\usecounter{l2}}}
\newcommand{\bdiamlist}{\begin{list}{$\diamond$}{}}
\newcommand{\bromalist}{\begin{list}{(\roman{l3})}{\usecounter{l3}}}

\newcommand{\prf}[1]{ \noindent {\em Proof:} \, #1 \hfill $\blacksquare$}

\newcommand{\argmin}{\mathop{\rm argmin}}
\newcommand{\argmax}{\mathop{\rm argmax}}
\newcommand{\diag}{\mathop{\mathrm{diag}}}
\newcommand{\tr}{\mathop{\rm Tr}}
\newcommand{\conv}{\mathop{\rm conv}}
\newcommand{\var}{\mathop{\rm Var}}
\renewcommand{\b}[1]{\ensuremath{\boldsymbol{\mathrm{#1}}}}
\newcommand{\Asur}{{\widetilde{A}}}
\newcommand{\Bsur}{{\widetilde{B}}}
\newcommand{\Gsur}{{\widetilde{L}}}
\newcommand{\Hsur}{{\widetilde{H}}}

\newcommand{\E}[1]{\b{\mu}_{{#1}}}
\newcommand{\Var}[1]{{\Sigma_{#1}}}

\newcommand{\bone}{\mathbf{1}}

\def\fic{v}
\def\contract{\Vcal}
\def\fictransform{Z}
\def\surrpol{\phi}
\def\surrpolset{\Phi}
\def\Xw{{P^w}}
\def\Xxi{{P^\xi}}
\def \couple{\Ccal}
\def \nodes{V}
\def \edges{E}
\def \gr{G}
\def \gnest{\gr_N}
\def \gcouple{\gr_C}
\def \enest{\edges_N}
\def \ecouple{\edges_C}
\def\support{\Wcal}
\def \pinfo{\zeta}
\def\proj{\Pi_\couple}
\def \gy{\gr}
\def \ey{\edges}
\def \gyset{\Gcal}
\def \lmi{\sf LMI}
\def \wtoX{L}

\def \add [#1]{\blue{#1}}
\def \remove [#1]{\red{#1}}
\def \replace [#1]#2{\red{#1} \blue{#2}}

\newcommand{\wl}[1]{\blue{(\textbf{Sam says:} #1)}}
\newcommand{\eb}[1]{\red{(\textbf{Eilyan says:} #1)}}

\newcommand{\rev}[1]{\black{#1}}
\newcommand{\conf}[1]{\black{#1}}

\graphicspath{{figures/}}

\usepackage{geometry}
\IEEEoverridecommandlockouts

\begin{document}

\title{ \LARGE \bf Decentralized Control of Constrained Linear Systems via Assume-Guarantee Contracts
}

\author{\vspace{.12in} Weixuan Lin$^{\dagger}$    \qquad Eilyan Bitar$^{\dagger}$  
\thanks{Supported in part by NSF grants ECCS-1351621 and IIP-1632124, and the Holland Phillips Trust.}
\thanks{$\dagger$ Weixuan Lin ({\tt\small wl476@cornell.edu}) and Eilyan Bitar ({\tt\small eyb5@cornell.edu}) are with the School of Electrical and Computer Engineering, Cornell University, Ithaca, NY, 14853, USA. }
}

\maketitle

\begin{abstract}
We consider the decentralized control  of a discrete-time, linear system subject  to exogenous disturbances and polyhedral constraints on the state and input trajectories. The underlying system is composed of a finite collection of  dynamically coupled subsystems, where each subsystem is assumed to have a dedicated local controller. The decentralization of information is expressed according to sparsity constraints on the state measurements that each local controller has access to. In this context, we investigate the design of decentralized controllers that are affinely parameterized in their measurement history. For problems with partially nested information structures, the optimization over such policy spaces is known to be convex.
Convexity is not, however, guaranteed under more general (nonclassical) information structures in which  the information available to one local controller can be affected by control actions that it cannot access or reconstruct.  
With the aim of alleviating the nonconvexity that arises in such problems, we propose an  approach to decentralized control design where the information-coupling  states  are effectively treated as disturbances whose trajectories are constrained to take values in ellipsoidal contract sets whose location, scale, and orientation are  jointly optimized with the underlying affine decentralized control policy. We establish a natural structural condition on the space of allowable contracts that facilitates the joint optimization over the control policy and the contract set via semidefinite programming.
\end{abstract}


\section{Introduction}

We investigate the design of affine decentralized control policies for 
stochastic discrete-time, linear systems that evolve over a finite horizon, and are subject to polyhedral constraints on the state and input trajectories. The computational tractability of such problems depends in part on their information structures  \cite{sandell1974solution, Tsitsiklis1985}. 
In particular, a decentralized control problem is said to have a \emph{nonclassical} information structure if the information available to one controller can be affected by the control actions of another that it cannot access or reconstruct. 
Under such information structures, the calculation of optimal decentralized control policies is known to be computationally intractable, because of the incentive for controllers to communicate with each other via the actions they undertake---the so called signalling incentive \cite{Witsenhausen1968, sandell1974solution, Tsitsiklis1985}.
To complicate matters further, there may be hard constraints  coupling the local actions and states of different controllers that must be jointly enforced
without explicit communication.
In this paper, we  address these challenges by relaxing the requirement that the decentralized controller be optimal with respect to the broad family of all causal policies, and instead search for suboptimal decentralized controllers that can be efficiently computed via convex programming methods.

\emph{Related Literature: \ }
There is a  related  literature  that leverages on techniques derived from tube-based model predictive control (MPC) to facilitate the design of  decentralized controllers for  constrained dynamical systems \cite{dunbar2007distributed, Dunbar2006, lucia2015contract, farina2012distributed, richards2004decentralized, richards2007robust, riverso2012tube, trodden2006robust, trodden2010distributed, hernandez2017distributed, trodden2014cooperative}.
Typically, these approaches rely on a decomposition of the decentralized control problem into a collection of decoupled local control problems by treating the coupling states and inputs associated with each subsystem's ``neighbors'' as independent exogenous disturbances that are assumed to take values in the given  state and input constraint sets.
Given the resulting collection of decoupled local control problems, centralized MPC methods can be applied to compute local control policies that are guaranteed to be feasible for each sub-problem.
Although decentralized control policies calculated according to such decomposition methods are guaranteed to be feasible for the full problem, they may result in behaviors that are overly conservative in terms of the  cost  they incur for a number of reasons.
\emph{First}, the treatment of the coupling states and inputs as independent disturbances ignores the potential dynamical coupling between these variables.
\emph{Second}, the over approximation of the coupling state and input trajectory sets by their corresponding state and input constraint sets will likely be very loose for many problem instances. More importantly, the over approximation of the coupling state and input trajectory sets in this manner ignores the fact that these sets depend on the control policy being used to regulate the system, and, therefore, neglects the possibility of co-optimizing their specification with the control policy.

\emph{Contribution: \ } 
We provide a computationally tractable method to calculate control policies that are guaranteed to be feasible for constrained decentralized control problems with nonclassical information structures.
Loosely speaking, the proposed approach seeks to neutralize the nonconvexity arising from the informational coupling between subsystems by treating the information-coupling states as  disturbances whose trajectories are ``assumed" to take values in a certain ``contract'' set. To ensure the satisfaction of this assumption, we impose a contractual constraint on the control policy that ``guarantees'' that the information-coupling states that it induces  indeed belong to said contract set.
Naturally, this approach yields an inner approximation of the original decentralized control design problem, where the conservatism of the resulting approximation depends on the specification of the contract set.
To limit the extent of the suboptimality that may result, we formulate a semi-infinite program to co-optimize the decentralized control policy with the location, scale, and orientation of an ellipsoidal contract set.
We establish a condition on the set of allowable contracts that facilitates the joint optimization of the control policy and the contract set via  semidefinite programming.

We note that there are several recent papers appearing in the literature that investigate a similar approach to decentralized control design via the co-optimization of control policies and contract sets \cite{trodden2017distributed, darivianakis2018decentralized}. Importantly, the techniques developed in these papers only permit the scaling and translation of a base contract set when co-optimizing it with the control policy. To the best of our knowledge, the method proposed in this paper provides the first systematic approach to co-optimize the control policy with the location, scale, and \emph{orientation} of the contract set,  expanding substantially the family of contracts that can be efficiently optimized over.

\emph{Notation: \ }   Let $\RR$ and $\RR_+$ denote the sets of real and non-negative real numbers, respectively.
Given a collection of vectors $x_1, \dots, x_N$, 
we let $(x_1, \dots, x_N)$ denote their vector concatenation in ascending order of their indices. Given an index set $J \subseteq \{1, \dots, N\}$, we let $x_J$ denote the  vector concatenation of the vectors $x_j$  for $j \in J$ in ascending order of their indices. Given a sequence $\{x(t)\}$ and time indices $s \leq t$, we let $x^{s:t} = (x(s), x(s+1), \dots, x(t) )$ denote its history from time $s$ to time $t$.
 Given a block matrix $A$, we let $[A]_{ij}$ denote its $(i,j)$-th block. We denote the trace of a square matrix $A$ by $\trace{A}$.
We denote the Minkowski sum of two sets $\Scal, \Tcal \subseteq \RR^n$ by $\Scal \oplus \Tcal := \{x + y \,  | \, x \in \Scal, \  y \in \Tcal \}$.

\section{Problem Formulation} \label{sec:formulation}

\subsection{System Model}

Consider a discrete-time, linear time-varying system consisting of $N$ dynamically coupled subsystems whose dynamics are described by
\begin{align} \label{eq:x_i}
x_i(t+1) = \sum_{j=1}^N  \left( A_{ij}(t) x_j(t)   + B_{ij} (t) u_j (t) \right) +  w_i (t),
\end{align}
for $i = 1,\dots, N$. 
We denote the \emph{local state}, \emph{local input}, and \emph{local disturbance} associated with  each subsystem $i$ at time $t$ by $x_i(t) \in \RR^{n_x^i}$, $u_i(t) \in \RR^{n_u^i}$, and $w_i (t) \in \RR^{n_x^i}$, respectively.
The system is assumed to evolve  over a finite time horizon $T$, and  the initial condition is assumed to be a random vector with known probability distribution.
In the sequel, we will work with a more compact representation of the full system dynamics given by
\begin{align*}
x(t+1) = A(t) x(t) + B(t) u (t) + w(t).
\end{align*}
Here, we denote by $x(t) := (x_1(t), .., x_N(t)) \in \RR^{n_x}$, $u(t) := (u_1(t), .., u_N(t)) \in \RR^{n_u}$, and $w (t) := (w_1(t), .., w_N(t)) \in \RR^{n_x}$ the full system state, input, and disturbance at time $t$. The dimensions of the system state and input are given by $n_x := \sum_{i=1}^N n_x^i$ and $n_u := \sum_{i=1}^N n_u^i$, respectively.

The input and disturbance trajectories are related to the state trajectory according to
\begin{align}
x = B u + \wtoX w, \label{eq:trajectory}
\end{align}
where $x$, $u$, and $w$ denote the system state, input, and disturbance trajectories, respectively.\footnote{The matrices $B$ and $\wtoX$ are specified in Appendix \ref{app:matrix}.}
They are defined by
\begin{align}
x & := (x(0), \dots, x(T))  \in \RR^{N_x}, \quad N_x := n_x(T+1), \label{eq:x_traj_def}\\ 
u & := (u(0), \dots, u(T-1))  \in \RR^{N_u}, \quad  N_u := n_u T,\label{eq:u_traj_def}\\
w & : = (w(-1), w(0), \dots, w(T-1))  \in \RR^{N_x} , \label{eq:w_traj_def}
\end{align}
where the initial component $w (-1)$ of the system disturbance trajectory is given by $w (-1 ) = x(0)$.
This notational convention will help simplify the specification of disturbance-feedback control policies in the sequel.

\subsection{Disturbance Model}

We model the disturbance trajectory $w$ as a zero-mean random vector whose support is an ellipsoid given by
\begin{align}
\support: = \left\{ z \in \RR^{N_x} \, \left|  \, z^\top \Sigma^{-1} z \leq 1\right. \right\}, \label{eq:ellipsoidal_support}
\end{align}
where the shape parameter $\Sigma \in \RR^{N_x \times N_x}$ is assumed to be symmetric and positive definite.
We let $M := \EE [w w^\top]$ denote the second moment matrix of the disturbance trajectory $w$. 
The matrix $M$ is guaranteed to be positive definite and finite-valued, as the support of $w$ is assumed to be an ellipsoid with a non-empty interior.

\subsection{System Constraints}

We consider a general family of polyhedral constraints on the state and input trajectories of the form 
\begin{align}
F_x x + F_u u + F_w w \leq g \qquad  \forall w \in \support, \label{eq:robust_con}
\end{align} 
where $F_x \in \RR^{m \times N_x}$, $F_u \in \RR^{m \times N_u}$, $F_w \in \RR^{m \times N_x}$, $g \in \RR^{m}$ are assumed to be given. Note that such constraints may couple states and inputs across subsystems and time periods.

\subsection{Information Structure}

We consider information structures that are specified according to \emph{sparsity constraints} on the state measurements that each controller has access to.
Specifically, we encode the pattern according to which information is shared between subsystems with a directed graph $\gr_I = (\nodes,\edges_I)$, which we refer to as the \emph{information graph} of the system. 
Here, the vertex set $\nodes = \{1, \dots , N \}$ assigns a distinct vertex $i$ to each subsystem $i$. Additionally, we include the directed edge
$(i,j) \in \edges_I$ if and only if  subsystem $j$ has access to subsystem $i$'s local state at each time $t$. 
We let $\nodes_I^-(i)$ denote the in-neighborhood of each subsystem $i \in \nodes$ in the information graph $\gr_I$.

Each subsystem is assumed to have access to the entire history of its local information up until and including time $t$.
More formally, we define the \emph{local information} available to each subsystem $i$ at time $t$ as
\begin{align} \label{eq:pol1}
z_i(t) : = \{ x_j^{0:t} \ | \ (j, i) \in \edges_I \}.
\end{align}
The local control input to each subsystem $i$ is restricted to be a causal function of its local information. That is, the local input to subsystem $i$ at time $t$ is of the form
\begin{align} \label{eq:pol2}
u_i(t)  = \gamma_i(z_i(t), t),
\end{align}
where $\gamma_i(\cdot, t)$ is a measurable function of the local information $z_i (t)$.
We  define the \emph{local control policy} for subsystem $i$ as $\gamma_i : = (\gamma_i(\cdot,0), \dots, \gamma_i(\cdot,T-1))$. We refer to the  collection of local control policies  $\gamma : = (\gamma_1, \dots, \gamma_N)$ as the \emph{decentralized control policy}, which relates the state trajectory $x$ to the input trajectory $u$ according to $u = \gamma (x)$.
Finally, we let $\Gamma$ denote the set of all decentralized control policies respecting the information constraints encoded in Eq. \eqref{eq:pol2}.

\subsection{Decentralized Control Design}

We consider the following family of constrained decentralized control design problems:
\begin{align}
\begin{alignedat}{8}
&\text{minimize} \qquad && \EE \left[ x^\top R_x x + u ^\top R_u u \right] \\
& \text{subject to} \qquad && \gamma \in \Gamma  \\
&&& \hspace{-.105in}\left. \begin{array}{l}
u = \gamma (x) \\
x = B u + \wtoX w  \\
F_x x + F_u u + F_w w \leq g 
\end{array}
\right\} \forall w \in \support.
\end{alignedat} \label{opt:decent}
\end{align}
Here, the cost matrices $R_x \in \RR^{N_x \times N_x}$ and $R_u \in \RR^{N_u \times N_u}$ are  assumed to be symmetric and positive semidefinite.
The tractability of the decentralized control design problem \eqref{opt:decent} depends  on the nature of the information structure.
In particular, if the information structure is \emph{partially nested}, then problem \eqref{opt:decent} can be equivalently reformulated (via the Youla parameterization) as a convex program in the space of disturbance feedback policies \cite{Ho1972}.
If, on the other hand, the information structure is nonclassical (i.e., not partially nested), then problem \eqref{opt:decent} is known to be computationally intractable, in general \cite{Tsitsiklis1985, Sandell1978, Mahajan2012_survey}.

\section{Information Decomposition} \label{sec:preliminaries}

The primary difficulty in solving decentralized control design problems stems from  the informational coupling that emerges when a subsystem's local information is affected by prior control actions that it cannot access or reconstruct. \rev{With the aim of isolating the effects of these actions on the  information available to each subsystem, we propose an   information  decomposition that partitions the local information available to each subsystem into a partially nested  subset (i.e., an information subset that is unaffected by control actions previously applied to the system) and its complement.}
This  information decomposition enables an equivalent reformulation of the decentralized control design problem where the control policy is expressed as an explicit function of the system disturbance and the so called \emph{information-coupling} states.
This reformulation will serve as the foundation for the contract-based approach to decentralized control design proposed in Section \ref{sec:feas_control}.

\subsection{Decomposition of Local Information}
We decompose the local information available to each subsystem according to a partition of  its in-neighbors in the   information graph $\gr_I$. More specifically, for each subsystem $i \in \nodes$, we let $$\Ncal(i) \subseteq \nodes_I^{-}(i)$$ denote the set of in-neighboring subsystems whose local state measurements contain information that is unaffected by the prior control actions of any subsystem. 
This requirement is satisfied if the local information of subsystem $i$ is such that it permits the reconstruction of all states and control actions directly affecting the local states of all subsystems belonging to $\Ncal(i)$. We denote the complement of this set by $\couple(i) : = \nodes_I^{-}(i) \setminus \Ncal(i)$ for each subsystem $i \in \nodes$.

With the goal of providing an explicit characterization of these sets, we first provide a characterization of the physical coupling between different subsystems as reflected by the block sparsity patterns of the system matrices $A$ and $B$. We describe this coupling in terms of a pair of directed graphs, $\gr_A := (\nodes, \edges_A)$ and $\gr_B : = (\nodes, \edges_B)$, whose edge sets  are defined according to
\begin{align*}
\edges_A := \{ (j,i) \in \nodes \times \nodes \ | \  \exists t = 0, \dots, T-1 \text{ s.t. } A_{ij} (t)\neq 0\},\\
\edges_B := \{ (j,i) \in \nodes \times \nodes \ | \  \exists t = 0, \dots, T-1 \text{ s.t. } B_{ij} (t)\neq 0\}.
\end{align*}
We let  $\nodes_A^- (i)$ and $\nodes_B^- (i)$ denote the in-neighborhoods associated with each node $i \in \nodes$ in  $\gr_A$ and $\gr_B$, respectively.

Building on these representations, we have the following definition that formalizes the class of information decompositions considered in this paper. For each subsystem $i \in \nodes$, define the set
\begin{align*}
\Ncal(i) := \{j \in  \nodes_I^- (i) \ | \  \text{\eqref{eq:cond1}, \eqref{eq:cond2} are satisfied}   \},
\end{align*}
where the above conditions  are given by
\begin{align} \label{eq:cond1}
\nodes_A^- (j) \subseteq \nodes_I^- (i) 
\end{align}
and
\begin{align} \label{eq:cond2}
\bigcup_{k \in \nodes_B^- (j)}  \nodes_I^- (k) \subseteq \nodes_I^- (i).
\end{align}
Condition \eqref{eq:cond1} requires that subsystem $i$ has access to all states that directly affect subsystem $j$'s state through the system dynamics. Condition \eqref{eq:cond2} requires that subsystem $i$ has access to the local information of each subsystem whose control actions directly affect subsystem $j$'s state. This  ensures that subsystem $i$ is able to reconstruct all control actions that directly affect subsystem $j$'s state.
Collectively, conditions \eqref{eq:cond1} and \eqref{eq:cond2} can be interpreted as a requirement on the \emph{local nesting of information}, in the sense that if $j \in \Ncal(i)$, then  subsystem $i$ is assumed to  have access to all states and control actions that directly affect subsystem $j$'s state through the state equation. As a result, subsystem $i$ can explicitly reconstruct the local disturbance $w_j (t-1)$ acting on any subsystem $j \in \Ncal(i)$ based only on its local information $z_i (t)$ as follows:
\begin{align*}
\nonumber w_j (t-1) = x_j (t) & - \sum_{k \in \nodes_A^- (j)} A_{jk}(t-1) x_k(t-1)   \\
& - \sum_{k \in \nodes_B^- (j)} B_{jk} (t-1) u_k (t-1) . 
\end{align*}

The local states of subsystems not belonging to $\Ncal(i)$, on the other hand, may contain information that can be influenced by prior control actions. We refer to these states as the \emph{information-coupling states} associated with subsystem $i$ at time $t$, denoting them by $x_{\couple(i)} (t)$ where $$\couple (i) := \nodes_I^- (i) \setminus \Ncal(i).$$ 
The collection of information-coupling states across all subsystems are denoted by the  $x_\couple (t) \in \RR^{n_x^\couple}$, where
\begin{align}
\couple := \bigcup_{i \in \nodes} \, \couple(i).
\end{align}
The trajectory of information-coupling states is denoted by
$$x_\couple := (x_\couple (0), \dots, x_\couple (T) ) \in \RR^{N_x^\couple},$$ 
where $N_x^\couple := n_x^\couple (T+1).$ 
Finally, it will be notationally convenient to express the mapping from the state trajectory $x$ to its subvector $x_\couple$ in terms of the projection operator  $\proj: \RR^{N_x} \to \RR^{N_x^\couple}$, where $x_\couple = \Pi_\couple x$.

\begin{remar}[Partially Nested Information] It can be shown that the given information structure is \emph{partially nested} if and only if the set of  information coupling states is empty, i.e., $\couple = \emptyset$. It is well known that such information structures permit the equivalent reformulation of problem \eqref{opt:decent} as a convex optimization problem in the space of disturbance-feedback control policies.
\end{remar}

\subsection{Control Input Reparameterization}

The proposed information decomposition suggests  a natural reparameterization of the control policy in terms of the following equivalent information set.

\begin{lemm}[Equivalence of Information] \label{lem:equi_info}
Define the information set $\zeta_i (t)$ according to
\begin{align*}
\zeta_i (t) := \{x_j^{0:t} | j \in \couple (i) \} \cup \{w_j^{-1: t-1}  | j \in \Ncal (i) \}.
\end{align*}
The sets $z_i (t)$ and $\zeta_i (t)$ are functions of each other for each subsystem $i$ and time $t$.
\end{lemm}

The proof of Lemma \ref{lem:equi_info} is omitted, as it mirrors that of \cite[Lemma 1]{Lin2016}.
Lemma \ref{lem:equi_info} suggests the following equivalent reparameterization of the local control input:
\begin{align}
u_i (t) = \surrpol_i (\pinfo_i (t), t), \label{eq:pol_equiv}
\end{align}
where $\surrpol_i (\cdot ,t)$ is a measurable function of its arguments. 
We let $\surrpol_i := ( \surrpol_i (\cdot, 0), \dots, \surrpol_i (\cdot, T-1) )$ and $\surrpol := (\surrpol_1, \dots, \surrpol_N)$ denote the reparameterized control policy associated with each subsystem $i \in \nodes$ and the full system, respectively.
With a slight abuse of notation, we express the input trajectory induced by  the reparameterized control policy $\surrpol$ as $$u = \surrpol (w, x_\couple).$$
Finally, we denote by $\surrpolset$ the set of reparameterized decentralized control policies that respect the information constraints encoded in Eq. \eqref{eq:pol_equiv}.

The reparameterization of the control input according to Eq. \eqref{eq:pol_equiv} results in the following  equivalent reformulation of the original decentralized control problem \eqref{opt:decent}:
\begin{align}
\begin{alignedat}{8}
&\text{minimize} \qquad && \EE \left[ x^\top R_x x + u ^\top R_u u \right] \\
& \text{subject to} \qquad && \surrpol \in \surrpolset  \\
&&& \hspace{-.105in}\left. \begin{array}{l}
u = \surrpol (w, x_\couple) \\
x = B u + \wtoX w  \\
F_x x + F_u u + F_w w \leq g 
\end{array}
\right\} \forall w \in \support.
\end{alignedat} \label{opt:decent_dFeedback}
\end{align}

Clearly, problem \eqref{opt:decent_dFeedback} remains nonconvex, in general, if the set of information-coupling subsystems is nonempty, i.e.,  $\couple \neq \emptyset$. In Section \ref{sec:feas_control}, we construct a convex inner approximation to problem \eqref{opt:decent_dFeedback}  where the information-coupling states are assumed to behave as disturbances with bounded support, and the control policy is constrained in a manner that ensures the consistency between  the assumed and actual behaviors of the information-coupling states.

\section{Decentralized Control Design via Contracts} \label{sec:feas_control}

In this section, we construct a convex inner approximation of the decentralized control design problem \eqref{opt:decent_dFeedback} via the introduction of an assume-guarantee contractual constraint on the information-coupling states $x_\couple$.
We do so by introducing a surrogate information structure in which the information-coupling states are modeled as fictitious disturbances that are  ``assumed" to take values in a  ``contract'' set.
To ``guarantee'' the satisfaction of this assumption, we impose a contractual constraint on the control policy requiring that the actual information-coupling states induced by the control policy belong to the contract set.
Given a fixed contract set, the resulting problem is a convex disturbance-feedback control design problem, whose feasible policies are guaranteed to be feasible for problem \eqref{opt:decent_dFeedback}.

\subsection{Surrogate Information}
We associate a \emph{fictitious disturbance}  $\fic_i (t) \in \RR^{n_x^i}$ with each subsystem $i \in \nodes$ and time $t = 0, \dots, T$. We let $\fic \in \RR^{N_x}$ denote the corresponding fictitious disturbance trajectory induced by these individual elements, which we model as a random vector whose support $\contract \subset \RR^{N_x}$ is assumed to be a convex and compact set.
We also assume that the  fictitious disturbance trajectory $v$ is independent of the system disturbance trajectory $w$.

Letting the collection of fictitious disturbances serve as surrogates for the information-coupling states, we define the \emph{surrogate local information} for subsystem $i$ as
\begin{align*}
\widetilde{\pinfo}_i (t) :=  \{v_j^{0:t} | j \in \couple (i) \} \cup \{w_j^{-1: t-1}  | j \in \Ncal (i) \}.
\end{align*}
Given a decentralized control policy $\surrpol \in \surrpolset$, the surrogate local information induces a \emph{surrogate control input} for each subsystem $i$ defined according to
\begin{align*}
\widetilde{u}_i (t) := \surrpol_i ( \widetilde{\pinfo}_i (t), t).
\end{align*}
Additionally, the \emph{surrogate input trajectory}  induced by the surrogate information structure is given by
\begin{align*}
\widetilde{u} := \surrpol (w, \fic_\couple),
\end{align*}
where $\fic_\couple := \proj \fic$.

\subsection{Surrogate Dynamics}

The treatment of the information coupling states as fictitious disturbances induces a \emph{surrogate system state} that evolves according to the following surrogate state equation:
\begin{align}
\nonumber \widetilde{x}_i (t+1) = & \sum_{j \in \nodes \setminus \couple(i) } A_{ij} (t) \widetilde{x}_j (t) + \sum_{j \in \couple  (i)} A_{ij} (t) \fic_j (t)  \\
&  + \sum_{j=1}^N B_{ij} (t) \widetilde{u}_j (t) + w_i (t), \label{eq:surr_state_i}
\end{align}
where $\widetilde{x}_i (t)$ denotes the surrogate state of subsystem $i$ at time $t$. We require that the initial condition of the surrogate system equal that of the true system, i.e.,  $\widetilde{x}_i (0) = x_i (0)$ for each subsystem $i$.
Moving forward, it will be convenient to express the  surrogate state dynamics in terms of trajectories as follows:
\begin{equation}
\widetilde{x} =  \Bsur \widetilde{u} +\Gsur w + \Hsur \fic_\couple,\label{eq:surrogate_state}
\end{equation}
where the matrices $\Bsur$, $\Gsur$, and $\Hsur$ are defined in Appendix \ref{app:matrix}.

We close this subsection with a lemma that establishes conditions  for the equivalence between the surrogate and actual state trajectories. We omit the proof, as it directly follows from the definition of the surrogate state equation \eqref{eq:surrogate_state}.

\begin{lemm} \label{lem:surrogate_state_eq}
\rev{Let $u  \in \RR^{N_u}$ and $w \in \RR^{N_x}$.} It holds that $x = Bu + \wtoX w$ if and only if $x = \Bsur u +\Gsur w + \Hsur x_\couple$.
\end{lemm}

\subsection{Assume-Guarantee Contracts} \label{sec:assume-guarantee}

Thus far, we have treated the information-coupling states as fictitious disturbances that are assumed to take values in a given  set $\contract_\couple$.
Leveraging on concepts grounded in assume-guarantee reasoning \cite{kwiatkowska2010assume, alur1999reactive},  we guarantee the satisfaction of this assumption by imposing a contractual constraint on the control policy, which  ensures that it induces information-coupling states that belong to $\contract_\couple$. We formalize the notion of an \emph{assume-guarantee contract} in the following definition.

\begin{definitio}[Assume-Guarantee Contract] \label{def:assume-guarantee}
A control policy $\surrpol \in \surrpolset$ is said to satisfy the \emph{assume-guarantee contract} speficied in terms of the \emph{contract set} $\contract_\couple \subseteq \RR^{N_x^\couple}$ if
\begin{align*}
\proj \widetilde{x} \in \contract_\couple \quad \forall \, (w , \fic_\couple  ) \in W \times \contract_\couple,
\end{align*}
where $ \widetilde{x} =  \Bsur \surrpol (w,   \fic_\couple) + \Gsur w + \Hsur  \fic_\couple $.
\end{definitio}
Here, the set $\contract_\couple$ is referred to as a \emph{contract set}, as it specifies the set that the information-coupling states are both assumed and required to belong to.
The satisfaction of the assume-guarantee contract guarantees that the \emph{surrogate} information-coupling states  $\widetilde{x}_\couple := \proj \widetilde{x}$ belong to the contract set.
In the following lemma, we show that the \emph{actual} information-coupling states that result under the policy $u = \surrpol (w, x_\couple)$ are guaranteed to belong to the contract set if the assume-guarantee contract is satisfied. 

\begin{lemm} \label{lem:assume-guarantee}
Let $\surrpol \in \surrpolset$ be a control policy that satisfies the assume-guarantee contract specified in terms of the contract set $\contract_\couple \subseteq \RR^{N_x^\couple}$. It follows that $\Pi_\couple x\in \contract_\couple$ for all $w \in \support$, where $x = B \surrpol (w, x_\couple) + \wtoX w$.
\end{lemm}

The proof of Lemma \ref{lem:assume-guarantee} is omitted due to space constraints.
In the following proposition, we provide an inner approximation of the decentralized control design problem \eqref{opt:decent_dFeedback} via the introduction of an assume-guarnatee contractual constraint. Its proof is omitted, as it follows directly from Lemma \ref{lem:assume-guarantee}.

\begin{propositio} \label{prop:feas}
Let $\surrpol \in \Phi$ be a feasible control policy for the following problem:
\begin{align}
\begin{alignedat}{8}
&\text{minimize} \qquad && \EE \left[ \widetilde{x}^\top R_x \widetilde{x} + \widetilde{u} ^\top R_u \widetilde{u} \right] \\
& \text{subject to} \qquad && \surrpol \in \surrpolset \\
&&& \hspace{-.105in}\left. \begin{array}{l}
\widetilde{u} = \surrpol (w,  \fic_\couple) \\
\proj \widetilde{x} \in \contract_\couple \\
\widetilde{x} =  \Bsur \widetilde{u} +\Gsur w + \Hsur \fic_\couple \\
F_x \widetilde{x} + F_u \widetilde{u} + F_w w \leq g  
\end{array}
\right\}  \forall (w, \fic_\couple) \in \support \times \contract_\couple,
\end{alignedat} \label{opt:contract}
\end{align}
It follows that $\phi$ is also  feasible for problem \eqref{opt:decent_dFeedback}.
\end{propositio}

Problem \eqref{opt:contract} is a convex disturbance feedback control design problem, given a fixed contract set $ \contract_\couple$. The choice of the contract set  does, however, play an important role in determining the performance of the control policies that it gives rise to. 
In Section \ref{sec:joint_opt}, we develop a systematic approach to enable the joint optimization of the contract set with the control policy via semidefinite programming.

\section{Policy-Contract Optimization} \label{sec:joint_opt}

In this section, we provide a semidefinite programming-based method to co-optimize the design of the decentralized control policy together with the contract set that constrains its design. 
As part of the proposed approach, we consider a restricted family of control policies that are affinely parameterized in both the disturbance and fictitious disturbance histories. We also parameterize the fictitious disturbance process as a causal affine function of a given (primitive) disturbance process---an approach that is similar in nature to the class of parameterizations that have been recently studied in the context of robust optimization with adjustable uncertainty sets \cite{zhang2017robust}. As one of our primary results in this section, we identify a structural condition on the family of allowable contract sets that permits the inner approximation of the resulting policy-contract optimization problem as a semidefinite program.

\subsection{Affine Control Policies}
We restrict our attention to affine decentralized disturbance-feedback control policies of the form
\begin{align}
\nonumber \widetilde{u}_i (t) =  u^o_i (t) &+ \sum_{ j \in \Ncal  (i) }   \sum_{s=-1}^{t-1} Q_{ij}^w (t,s+1) w_j (s)  \\
&+ \sum_{ j \in \couple (i)} \sum_{s=0}^t  Q_{ij}^\fic (t,s) \fic_j (s), \label{eq:aff1}
\end{align}
for $t = 0, \dots, T-1$ and $i = 1, \dots, N$. Here, $u^o_i (t)$ denotes the open-loop control input, and the matrices $Q_{ij}^w (t,s+1) $ and $ Q_{ij}^\fic (t,s)$ denote the feedback control gains.
The affine control policy specified in Eq. \eqref{eq:aff1} can be expressed in terms of trajectories as
\begin{align}
\widetilde{u} = u^o +  Q^w  w + Q^\fic \fic, \label{eq:aff2}
\end{align}
where the gain matrices $Q^w$ and $Q^\fic$ are both $T \times (T+ 1)$ block matrices, whose $(t, s)$-th block is defined according to
\begin{align}
[Q^w (t,s) ]_{ij} &= \begin{cases}
Q_{ij}^w (t,s) & \text{if } j \in \Ncal (i) \ \text{and} \ t \geq s, \\
0 & \text{otherwise},
\end{cases}  \label{eq:Qw}\\
[Q^\fic (t,s)] _{ij} &= \begin{cases}
Q_{ij}^\fic (t,s) & \text{if } j \in \couple (i)  \ \text{and} \ t \geq s, \\
0 & \text{otherwise}.
\end{cases}  \label{eq:Qfic}
\end{align}
for $i, j = 1, \dots, N$.
We let $\Qcal_N$ and $\Qcal_C$ denote the matrix subspaces  respecting the block sparsity patterns specified according to Eqs. \eqref{eq:Qw} and \eqref{eq:Qfic}, respectively.

\subsection{Affine Parameterization of the Fictitious Disturbance}

We focus our analysis on fictitious disturbances that are expressed according to affine transformations of a \emph{primitive disturbance}. Such a parameterization yields contract sets that have adjustable location, scale, and orientation.
Specifically, we let the random vector $\xi$ denote the \emph{primitive disturbance trajectory}, which is assumed  to be  an i.i.d.  copy of the system disturbance trajectory $w$. We parameterize the fictitious disturbance trajectory affinely in the primitive disturbance as
\begin{align}
\fic := \overline{\fic} + \fictransform \xi. \label{eq:zeta}
\end{align}
Here, the parameters  $\overline{\fic} \in \RR^{N_x}$ and $\fictransform \in \RR^{N_x \times N_x}$ can be adjusted to control the shape of the resulting contract set $\Vcal_\couple$, which takes the form of
\begin{align}
\contract_\couple = \proj \left( \overline{\fic} \oplus \fictransform \support \right). \label{eq:contract}
\end{align}
Throughout the paper, we will restrict our attention to transformations \eqref{eq:zeta} in which the matrix parameter $\fictransform$ is both lower triangular and invertible. We denote the set of all such matrices by  $\Zcal \subset \RR^{N_x \times N_x}$.

The specification of the fictitious disturbance according to Eq. \eqref{eq:zeta} induces the following the surrogate control input:
\begin{align}
\widetilde{u} = u^o + Q^\fic \overline{\fic} +  Q^w  w + Q^\fic \fictransform \xi. \label{eq:aff3}
\end{align}
We eliminate the bilinear terms in Eq. \eqref{eq:aff3} through the following  the change of variables:
\begin{align}
 \overline{u}:= u^o + Q^\fic \overline{\fic}   \ \text{ and } \ Q^\xi := Q^\fic \fictransform. \label{eq:change_of_var}
\end{align}
This change of variables gives rise to a reparameterization of the surrogate input trajectory as
\begin{align}
\widetilde{u} = \overline{u} +  Q^w  w + Q^\xi \xi, \label{eq:aff4}
\end{align}
where the matrix $Q^\xi \in \RR^{N_u \times N_x}$ must satisfy the sparsity constraint 
$$Q^\xi \fictransform^{-1} \in \Qcal_C$$
in order to ensure the satisfaction of the original sparsity constraint that $Q^\fic \in \Qcal_C$.

The parameterization of the contract set and control policy in this manner permits their co-optimization as follows:
\begin{align}
\begin{alignedat}{8}
&\text{minimize} \qquad && \EE \left[ \widetilde{x}^\top R_x \widetilde{x} + \widetilde{u} ^\top R_u \widetilde{u} \right] \\
& \text{subject to} \qquad &&Q^w \in \Qcal_N,  Q^\xi\in \RR^{N_u \times N_x} ,  \fictransform \in \Zcal \\
&&& \overline{u} \in \RR^{N_u}, \  \overline{\fic} \in \RR^{N_x} , \\
&&& Q^\xi \fictransform^{-1} \in \Qcal_C \\
&&& \hspace{-.105in} \left. \begin{array}{l}
\fic  =  \overline{\fic}  + \fictransform \xi\\
\widetilde{u} = \overline{u} +  Q^w  w + Q^\xi \xi \\
\widetilde{x} =  \Bsur \widetilde{u} +\Gsur w + \Hsur \fic_\couple  \\
\proj \widetilde{x} \in \proj  \left( \overline{\fic} \oplus \fictransform \support \right) \\
F_x \widetilde{x} + F_u \widetilde{u} + F_w w \leq g
\end{array}
\right\} \forall (w, \xi) \in \support^2,
\end{alignedat} \label{opt:inner0_joint}
\end{align}
where $\support^2 := \support \times \support$. 
Problem \eqref{opt:inner0_joint} is a nonconvex semi-infinite program, where the nonconvexity is due to  the sparsity constraint on the matrix $Q^\xi \fictransform^{-1}$ and the contractual constraint on the affine control policy. In what follows, we provide convex inner approximations of these constraints, which  yield an inner approximation of problem \eqref{opt:inner0_joint} as a semidefinite program.

\subsection{Restricting the Contract Set}

In what follows, we introduce an additional restriction  on the set of allowable matrix parameters $\fictransform$ that guarantees the invariance of the subspace $\Qcal_C$ under multiplication by such  matrices. This permits the equivalent reformulation of the bilinear constraint $Q^\xi \fictransform^{-1} \in \Qcal_C$ as $Q^\xi \in \Qcal_C$.

Specifically, we require that the matrix $\fictransform$ be of the form
\begin{align}
\fictransform = \lambda I - Y, \label{eq:fictransform}
\end{align}
where $\lambda \geq 1$ is scalar parameter and $Y \in \RR^{N_x \times N_x}$ is a $(T+1) \times (T+1)$ strictly block lower triangular matrix  of the form
\begin{align}
Y = \bmat{0 \\ Y(1,0) & 0 \\ 
\vdots &   \ddots  & \ddots\\
Y (T, 0 ) &  \cdots & Y(T, T-1) & 0}. \label{eq:Y}
\end{align}
Furthermore, each block of the matrix $Y$ is an $N \times N$ block matrix, whose $(i, j)$-th block is of dimension $n_x^i \times n_x^j$.   We impose an additional restriction on the structure of the matrix $Y$ in the form of sparsity constraints (that reflect the pattern of informational coupling between subsystems) on each of its blocks.

More specifically, we encode the pattern of informational coupling between subsystems according to a directed graph $\gcouple := (\nodes, \ecouple)$, whose  directed edge set $\ecouple$ is defined as
\begin{align*}
\ecouple := \{(j, i) \in \edges_I \, | \, j \in \couple (i) \}.
\end{align*}
We let $\nodes_C^+ (i)$ denote the out-neighborhood of a node $i \in \nodes $ in the graph $\gcouple$.
Using this graph, we impose a sparsity constraint  on each  block of the matrix $Y$ of the form:
\begin{align}
[Y (t, s) ]_{ij} = 0 \quad \text{if } \ \nodes_C^+ (i) \nsubseteq \nodes_C^+ (j)  \label{eq:subspace_Y}
\end{align}
for all $i,j = 1, \dots, N$, and $t, s = 0,\dots, T$.
We  let $\Ycal (\gcouple)$ denote the subspace of all matrices that respect the sparsity constraints implied by Eqs. \eqref{eq:Y} and \eqref{eq:subspace_Y}.

We have the following result establishing the invariance of the subspace $\Qcal_C$ under multiplication by matrices 
$Y \in \Ycal (\gcouple)$.

\begin{lemm} \label{lem:QY} If $Q \in \Qcal_C$ and $Y \in \Ycal (\gcouple)$, then $Q Y \in \Qcal_C$
\end{lemm}
\prf{
The sparsity constraint $QY \in \Qcal_C$ is satisfied if the matrix $Q (t, s) Y (s, r)$ satisfies the sparsity constraint
\begin{align*}
[Q (t, s) Y (s, r)]_{ij} = 0 \quad \forall i \notin \nodes_C^+ (j)
\end{align*}
for all times $r, \; s, \; t$ satisfying $0 \leq r < s \leq t \leq T-1$.
We prove this claim by showing that $[Q (t, s) Y (s, r)]_{ij} \neq 0 $ implies $ i \in \nodes_C^+ (j)$. The condition that $[Q (t, s) Y (s, r)]_{ij} \neq 0$  implies that there exists $k \in \nodes$ such that the blocks $[Q (t, s)]_{ik}$ and $[Y (s, r ) ]_{kj}$ are both nonzero. 
The fact that $[Q (t, s)]_{ik}$ is nonzero implies that $i \in \nodes_C^+ (k)$, as the matrix $Q$ satisfies $Q \in \Qcal_C$. The fact that $[Y (s, r ) ]_{kj}$ is nonzero implies that $\nodes_C^+ (k) \subseteq \nodes_C^+ (j)$, as the matrix $Y$ satisfies $Y \in \Ycal (\gcouple)$. The desired result follows. 
}

We have the following result as an immediate consequence of Lemma \ref{lem:QY}.

\begin{lemm} \label{lem:change_of_variable}
Let  $Y \in \Ycal (\gcouple)$  and $\lambda \in [1, \infty)$. It follows that
\begin{align*}
\Qcal_C = \left\{ Q^\xi ( \lambda I - Y)^{-1} \, | \, Q^\xi \in \Qcal_C  \right\}.
\end{align*}
\end{lemm}
It follows from Lemma \ref{lem:change_of_variable} that the constraint $Q^\xi Z^{-1} \in \Qcal_C$ is equivalent to $Q^\xi \in \Qcal_C$ 
if  $Z = \lambda I - Y$, where $Y \in \Ycal (\gcouple)$  and $\lambda \geq 1$.

\subsection{Semidefinite Programming Approximation}

To lighten notation, we write the surrogate state trajectory $\widetilde{x}$ more compactly as
\begin{align*}
\widetilde{x} = \overline{x} +  \Xw  w +  \Xxi \xi ,
\end{align*}
where $\overline{x}  := \Bsur \overline{u} + \Hsur \proj \overline{\fic}$, $\Xw:= \Bsur Q^w + \Gsur$,  and $\Xxi:= \Bsur Q^\xi + \Hsur \proj (\lambda I - Y)$. 

We first address  the robust linear constraints in problem \eqref{opt:inner0_joint}. The following result provides an equivalent reformulation   as second-order cone constraints. Its proof is omitted, as it is an immediate consequence of the identity  $\sup_{w \in \support} c^\top w = \lVert \Sigma^{1/2} c \rVert_2$ for all $c \in \RR^{N_x}$.
\begin{lemm} \label{lem:robust_linear}
The semi-infinite constraint $F_x \widetilde{x} + F_u \widetilde{u} + F_w w \leq g $ for all $(w,  \xi) \in \support^2 $ is satisfied if and only if
\begin{align}
\nonumber &  \left\lVert \Sigma^{1/2} e_i^\top  ( F_x \Xw + F_u Q^w + F_w ) \right\rVert_2 \\
\nonumber &\qquad + \left\lVert \Sigma^{1/2} e_i^\top ( F_x \Xxi + F_u Q^{\xi}) \right\rVert_2 \\
& \qquad \leq e_i^\top ( g - F_x \overline{x} - F_u \overline{u} ), \quad i = 1, \dots, m, \label{eq:robust_linear}
\end{align}
where $e_i$ is the $i$th standard basis vector in $\RR^m$. 
\end{lemm}

We now address the nonconvexity that stems from the
contractual constraint in problem \eqref{opt:inner0_joint}. First, notice  that the contractual constraint is equivalent to the following set containment constraint
\begin{align}
\proj \left( \overline{x} \oplus \Xw \support \oplus  \Xxi \support \right)    \subseteq \proj  \left( \overline{\fic} \oplus  \fictransform \support  \right). \label{eq:containment}
\end{align}
The set containment constraint \eqref{eq:containment} amounts to requiring that the Minkowski sum of two ellipsoids be contained within another ellipsoid.  It follows from \cite{durieu2001multi}[Theorem 4.2] that this class of set containment constraints can be approximated from within by a quadratic matrix inequality. Through an application of Schur's Lemma, one can approximate the resulting quadratic matrix inequality from within by a linear matrix inequality. We summarize  the resulting inner approximation in the following lemma.

\begin{lemm} \label{lem:matrix_inequality}
The set containment constraint \eqref{eq:containment} is satisfied if there exists a scalar $\beta \in [0, \lambda]$ such that
\begin{align}
&\proj \left(   \overline{x} - \overline{\fic} \right) =0 , \label{eq:LMI0}\\
&\bmat{\proj \widetilde{\Sigma} \proj^\top & \proj \Xw & \proj \Xxi \\
\Xw^\top \proj^\top & \beta \Sigma^{-1} & 0 \\
\Xxi ^\top \proj^\top & 0 & (\lambda - \beta) \Sigma^{-1} } \succeq 0, \label{eq:LMI}
\end{align}
where  $\widetilde{\Sigma} = \lambda \Sigma - Y \Sigma - \Sigma Y^\top.$
\end{lemm}

By applying Lemmas \ref{lem:change_of_variable}--\ref{lem:matrix_inequality}, one can approximate the nonconvex semi-infinite program \eqref{opt:inner0_joint} from within as the following finite-dimensional semidefinite program.
\begin{propositio}
Each feasible solution to the following semidefinite program is feasible for problem \eqref{opt:inner0_joint}:
\begin{align}
\begin{alignedat}{8}
&\text{minimize} \qquad && {\rm Tr} \left( \Xxi^\top R_x  \Xxi M + \Xw^\top R_x  \Xw M \right) \\
&&& + {\rm Tr} \left({Q^w}^\top R_u Q^w M +  {Q^\xi}^\top R_u Q^\xi M  \right) \\
&&&  + \overline{x}^\top R_x  \overline{x} + \overline{u}^\top R_u  \overline{u} \\[3pt]
& \text{subject to} \qquad &&Q^w \in \Qcal_N, Q^\xi \in \Qcal_C , Y \in \Ycal (\gcouple) , \\
&&& \overline{u} \in \RR^{N_u}, \overline{\fic}, \ \overline{x} \in \RR^{N_x} , \lambda, \ \beta \in \RR_+ ,\\
&&&  \Xw, \ \Xxi  \in \RR^{N_x \times N_x}, \\
&&& \lambda \geq \max \{ 1, \beta \}, \\
&&& \overline{x} =  \Bsur \overline{u} + \Hsur \proj \overline{\fic} \\
&&& \Xw = \Bsur Q^w + \Gsur \\
&&& \Xxi = \Bsur Q^\xi + \Hsur \proj (\lambda I - Y) \\
&&& \eqref{eq:robust_linear}, \ \eqref{eq:LMI0}, \ \eqref{eq:LMI}.
\end{alignedat} \label{opt:inner4_joint}
\end{align}
\end{propositio}

The decision variables for problem \eqref{opt:inner4_joint} are the matrices $Q^w$, $Q^\xi$, $Y$, $\Xw$, $\Xxi$, the vectors $\overline{u}$, $\overline{v}$, $\overline{x}$, and the scalars $\lambda$ and $\beta$. 
Problem \eqref{opt:inner4_joint} is a convex inner approximation of the reformulated decentralized control design problem \eqref{opt:decent_dFeedback}, in the sense that each feasible solution of problem \eqref{opt:inner4_joint} can be mapped to a feasible affine control policy for problem \eqref{opt:decent_dFeedback} via the change of variables specified in  \eqref{eq:change_of_var}.
The decentralized control policies that this approximation gives rise to are suboptimal, in general.  Bounds on their suboptimality, however, can be efficiently calculated using information-based convex relaxations \cite{lin2019convex}.

\section{Conclusion} \label{sec:conclusions}

We provide a method to compute feasible control policies for constrained decentralized control design problems by leveraging on the concept of assume-guarantee contracts. At the heart of this approximation is the treatment of information-coupling states as fictitious disturbances that are ``assumed" to take values in a contract set. We ``guarantee" the inclusion of the information-coupling states in the contract set by imposing an assume-guarantee contractual constraint on the control policy.
The introduction of such assume-guarantee contracts gives rise to an inner approximation of the decentralized control design problem, whose quality depends on the specification of the contract set. We provide a method of co-optimizing the decentralized control policy with the location, scale, and orientation of the contract set via semidefinite programming.


\bibliographystyle{IEEEtran}
\bibliography{decentralized_bib}{\markboth{References}{References}}

\begin{appendices}

\section{Matrix Definitions} \label{app:matrix}

\rev{
Define the matrices $\Asur (t)$ and $ \Hsur (t)$ according to
\begin{align*}
\Asur_{ij} (t) &= \begin{cases}
A_{ij} (t)  & \text{if } j \in \nodes \setminus \couple (i) ,\\
0 & \text{otherwise,}
\end{cases} \\
\Hsur_{ij} (t) &= A_{ij} (t) - \Asur_{ij} (t),
\end{align*}
where $i, j \in \nodes$.}
The matrices $(B, \wtoX)$ in Eq. \eqref{eq:trajectory} and the matrices $(\Bsur, \Gsur)$ in Eq. \eqref{eq:surrogate_state} are defined according to
\begin{align*}
 B & := \bmat{0 \\ A_1^1 B(0) &0 \\ A_1^2 B(0) & A_2^2 B(1) & 0 \\  \vdots &  & \ddots & \ddots \\ \vdots &&&\ddots &0 \\ A_1^T B(0) & A_2^T B(1)  & \cdots & \cdots &A_T^T B(T-1)},\\
 \Bsur & := \bmat{0 \\ \Asur_1^1 B(0) &0 \\ \Asur_1^2 B(0) & \Asur_2^2 B(1) & 0 \\  \vdots &  & \ddots & \ddots \\ \vdots &&&\ddots &0 \\ \Asur_1^T B(0) & \Asur_2^T B(1)  & \cdots & \cdots &\Asur_T^T B(T-1)},\\
 \wtoX & := \bmat{A_0^0  \\ A_0^1 & A_1^1  \\  \vdots &   & \ddots \\  A_0^T & A_1^T   &    \cdots &A_T^T }, \
  \Gsur  := \bmat{\Asur_0^0  \\ \Asur_0^1 & \Asur_1^1  \\  \vdots &   & \ddots \\  \Asur_0^T & \Asur_1^T   &    \cdots &\Asur_T^T }, 
\end{align*}
where $A_s^t := \prod_{r=s}^{t-1} A(r)$ and $\Asur_s^t := \prod_{r=s}^{t-1} \Asur(r)$ for $s < t$, and $A_t^t = \Asur_t^t = I $. Additionally, the matrix $\Hsur$ in Eq. \eqref{eq:surrogate_state} is defined as $\Hsur := H \proj ^\top$, where
\[   H  := \bmat{0 \\ \Asur_1^1 \Hsur(0) &0 \\ \Asur_1^2 \Hsur(0) & \Asur_2^2 \Hsur(1) & 0 \\  \vdots &  & \ddots & \ddots \\ \vdots &&&\ddots &0 \\ \Asur_1^T \Hsur(0) & \Asur_2^T \Hsur(1)  & \cdots & \cdots &\Asur_T^T \Hsur(T-1) & 0}. \]

\end{appendices}

\end{document}